# Complete graphs: the space of simplicial cones, and their path tree representation

*Amos Ron*\*, *Shengnan (Sarah) Wang*\*\*

ABSTRACT

Let $G$ be a complete graph with $n+1$ vertices. In a recent paper of the authors, it is shown that the path trees of the graph play a special role in the structure of the truncated powers and partition functions that are associated with the graph. Motivated by the above, we take here a closer look at the geometry of the simplicial cones associated with the graph, and the role played by those simplicial cones that are generated by path trees. It is shown that the latter form a basis for the linear space spanned by the former, and that the representation of a general simplicial cone by path tree cones can be deduced by examining partial orders induced by rooted trees. While the problem itself is geometrical and its solution is combinatorial, the proofs rest with multivariate spline theory.

**Authors' affiliation and address:**

\* Computer Sciences Department  
University of Wisconsin-Madison  
1210 West Dayton Street  
Madison WI 53706  
`amos@cs.wisc.edu`

\*\* Morgridge Institute for Research  
Madison WI 53715  
`swang@morgridge.org`  
*and*  
Department of Electrical & Computer Engineering  
University of Illinois - Urbana Champaign  
Urbana IL 61801

This work was supported by the National Science Foundation under Grant DMS-1419103.

## 1. Introduction

It was recently shown in [8] that the path trees in an acyclic graph play a special role in the decomposition of the associated truncated powers and partition functions into 'cone polynomials', i.e., summands that are, each, a single polynomial supported on a single cone. Motivated by this decomposition, we attempt in this paper to gain an insight into the role played by path trees in the geometry of the simplicial cones associated with a complete graph. We have chosen the complete graph setup since then the *entire* geometry is captured by the path trees. Our setup consists of a space of piecewise-constants, defined on the geometry of the above simplicial cones, and the goal is to get a complete grasp of the structure of that space. It turns out that this structure hinges on intricate relationships among the spanning trees of the graph, and that, indeed, the understanding of those relationships is best achieved by utilizing the path trees of the graph.

A 'take home' summary of the results of this paper may go as follows. Relationships among the simplicial cones associated with a graph hinge on the different partial orders on the vertices that the spanning trees of the graph induce. The path trees are unique here, since they are the only ones to induce full order.

Let $G$ be an acyclic complete graph with vertex set $[0{:}n]$. With

$$(e_i)_{i=1}^n$$

the standard basis for $\mathbb{R}^n$, and with $e_0 := 0 \in \mathbb{R}^n$, we assume that the orientation is given by

$$i \to j \iff i < j, \quad i,j \in [0{:}n].$$

That means that the vector representation of every edge $(i,j) \in G$ is

$$x = e_j - e_i, \quad 0 \leq i < j \leq n.$$

We identify in what follows the edges with their vector representation, and the graph itself with the set of vectors that represent the edges. Other orientations yield results that are completely analogous to the ones below. Let

$$\mathbb{B}(G)$$

be the set of spanning trees of $G$, oriented as in $G$. Recall that

$$\#\mathbb{B}(G) = (n+1)^{n-1}.$$

Given $T \in \mathbb{B}(G)$, the positive hull

$$\mathrm{pos}(T) \subset \mathbb{R}^n$$

is a simplicial cone whose support function is denoted by

$$\chi_T.$$

Our interest is in the space
$$\mathcal{T}(G) := \operatorname{span}\{\chi_T : T \in \mathbb{B}(G)\}.$$

It is a piecewise constant space that encodes the geometry of the underlying truncated powers and partition functions: whatever multiplicity one assigns to the edges of the complete graph, the domains of polynomiality of the truncated power coincide with the domains of constant values of a generic function in $\mathcal{T}(G)$. Since the geometry (in contrast with the algebra) as investigated here does not depend on the multiplicities of the edges, we assume throughout this paper that each edge is simple, i.e., appears with multiplicity 1.

Our first result in this investigation is

**Proposition 1.1.**
$$\dim(\mathcal{T}(G)) = n!.$$

□

In addition to the cone $\operatorname{pos}(T)$ that is defined with respect to the $G$-orientation of $T$, the spanning tree $T$ also induces a partial ordering on $[0{:}n]$ as follows. One first re-orients $T$ as a tree rooted at 0. We denote this orientation of $T$ as
$$T_{on}$$
and refer to it as the *natural orientation of $T$*. The partial ordering $\prec_T$ is then induced by $T_{on}$: $i \prec_T j$ iff $i \neq j$ and the path in $T$ from 0 to $j$ contains $i$. A *path tree* $T$ is a spanning tree whose natural orientation is
$$0 \to s(1) \to \ldots \to s(n),$$
with $s \in S_n$, i.e., $s$ is a permutation of the vertex set $[1{:}n]$. We write $T_s$ for the above path tree when oriented as in $G$, and $T_{s,on}$, thus, for the same tree oriented as above. We then denote
$$\chi_s := \chi_{T_s}, \quad s \in S_n.$$

We have:

**Theorem 1.2.** *The path simplicial cone functions*
$$\{\chi_s : s \in S_n\}$$
*form a basis for $\mathcal{T}(G)$.*

Thus, every $\chi_T$, $T \in \mathbb{B}(G)$, can be written uniquely as
$$\chi_T = \sum_{s \in S_n} c(T,s) \chi_s,$$
for some real coefficients $c(T,s)$. We have

**Proposition 1.3.** *For each $T \in \mathbb{B}(G)$ and for each $s \in S_n$,*
$$c(T,s) \in \{0, \pm 1\}.$$

□

Next, every $s \in S_n$ induces a new orientation on $G$ that we denote as $G^s$:

$$x = e_{s(j)} - e_{s(i)} \in G^s \iff i < j.$$

Here, $s(0) := 0$, so 0 is always the *source* of $G^s$, i.e., it does not have *inflow* edges. Note that every complete, acyclic, graph has exactly one source. The re-orientation of any $Y \subset G$ according to the orientation $G^s$ is denoted by

$$Y^s.$$

Thus, in this language, $T_s^s = T_{s,on}$: the natural orientation of $T_s$ is its orientation in $G^s$.

Given $s \in S_n$, and $T \in \mathbb{B}(G)$, the *distortion*

$$\mathbf{d}(T)$$

is the number of edges in $T$ whose orientations in $T$ and $T_{on}$ differ. Finally, we say that two trees, $T, T' \in \mathbb{B}(G)$ are *compatible* if their partial orderings $\prec_T$ and $\prec_{T'}$ are compatible, i.e., if there exists a full ordering of $[0{:}n]$ that is compatible with both, and define

$$\mathbb{B}(G) \times \mathbb{B}(G) \ni (T, T') \mapsto C(T, T') := \begin{cases} 1, & T, T' \text{ are compatible,} \\ 0, & \text{otherwise.} \end{cases}$$

Our final result in this paper determines the exact value of each $c(T, s)$:

**Theorem 1.4.** *Let $T$ and $s$ be as above. Then:*

$$c(T, s) = (-)^{\mathbf{d}(T) + \mathbf{d}(T_s)} C(T, T_s).$$

Note that, according to this theorem, the cones of two permutations $s, s' \in S_n$ of different distortion parity are never summed up in the representation of any $\chi_T$.

**Example.** Let $T = (e_i : i \in [1 : n])$. Then $T$ is compatible with any other spanning tree, and $\mathbf{d}(T) = 0$, hence

$$\chi_T = \sum_{s \in S_n} (-)^{\mathbf{d}(T_s)} \chi_s.$$

For any other tree $T$, the cone function $(-)^{\mathbf{d}(T)} \chi_T$ is obtained by dropping suitable terms from the above sum.

**Comment.** While the distortion $\mathbf{d}(T_s)$ does not depend on the simplicial cone $\chi_T$ we represent, it does depend on the orientation of $G$. Thus, for a fixed orientation of $G$, any two path simplicial cones have always the same sign relationship in their representation of the $\chi_T$ function. For example, two path cones with different parity of distortion will always appear with opposite signs of the coefficients, provided that both coefficients are not 0. However, if we change the orientation of $G$, the same path cones may now have the same parity of distortion, hence their signs will now agree in each representation.

## 2. Examples

We represent $s \in S_n$ as a sequence of integers. For example, the sequence [4123], corresponds to $s(1) = 4$, $s(2) = 1$, etc., hence to the path tree

$$T_{s,on} : 0 \to 4 \to 1 \to 2 \to 3.$$

To recall, we always extend $s$ with $s(0) = 0$, hence [4123] and [04123] represent the same $s \in S_n$, and the same path tree $T_s$.

### 2.1. $\mathcal{T}(G)$ for $n = 2$.

Figure 1 shows the geometry of this case. There are two path cones, marked green and red. The remaining cone $\chi$ (blue) corresponds to $T = \{e_1, e_2\}$. Obviously,

$$\chi = \chi_{[12]} - \chi_{[21]}.$$

This is a special case of the example at the end of the previous section. Indeed, for the orientation $s = [12]$, we have that $T_s = T_{s,on}$, hence $\mathbf{d}(T_{[12]}) = 0$. As for the orientation $[21]$, the edge $e_2 - e_1$ is a mismatch when comparing $T_s$ and $T_{s,on}$, hence $\mathbf{d}(T_{[21]}) = 1$.

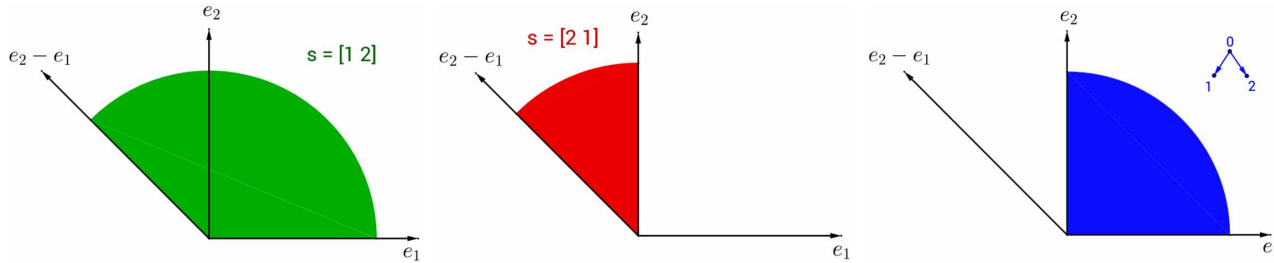

**Figure 1**: The 2D case

### 2.2. $\mathcal{T}(G)$ for $n = 3$.

There are $4^2 = 16$ members in $\mathbb{B}(G)$, when $n = 3$. Let us first examine the partial orders induce by those trees. We naturally skip the six path trees: they induce a full order, which is compatible only with the full order of themselves. So, Theorem 1.4 states correctly then the trivial assertion that the cone function $\chi_s$, $s \in S_n$, is a linear combination of itself.

In the illustrations below, we show a cross-cut of the 3D $pos(G)$, and use the six edges to denote the intersection points of their rays with the planar cut.

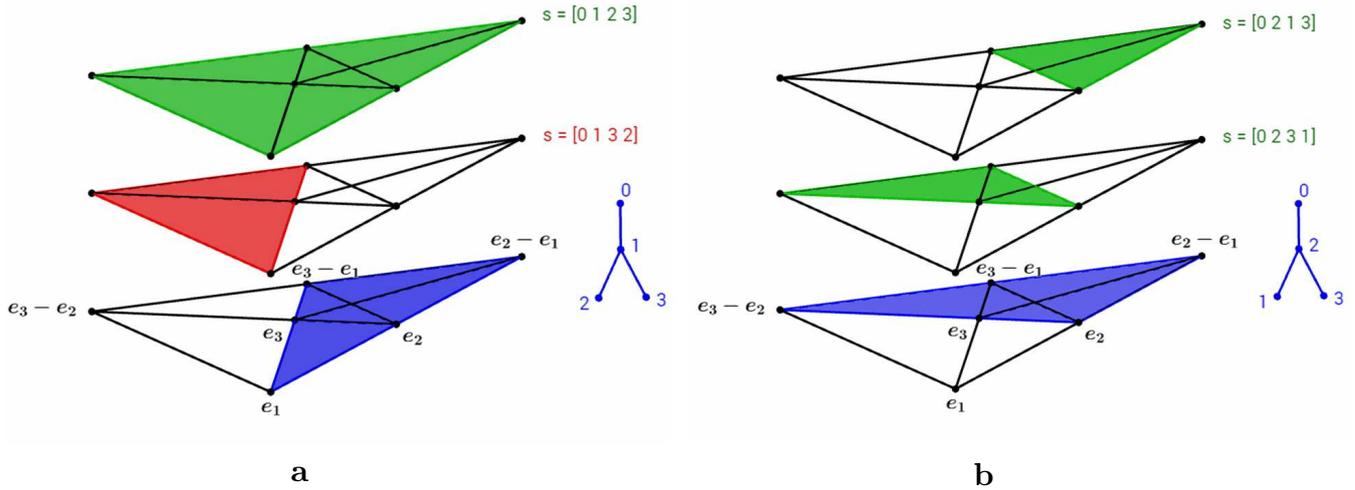

**Figure 2**: Representations by two path cones

Of the ten non-path cones, three are of the form shown in Figure 2. The partial order can then be completed to a full one in two different ways. For the specific tree $T$ shown in Figure 2(a), $\mathbf{d}(T) = 0$, and one needs just to know how to orient the edge $e_3 - e_2$ in order to obtain a full order. For the order $s = [0123]$, we have $\mathbf{d}(T_s) = 0$, and for $s' = [0132]$, we have $\mathbf{d}(T_{s'}) = 1$. The induced decomposition is shown in Figure 2(a): the positive cone in green, the negative in red, and the cone $\chi_T$ in blue.

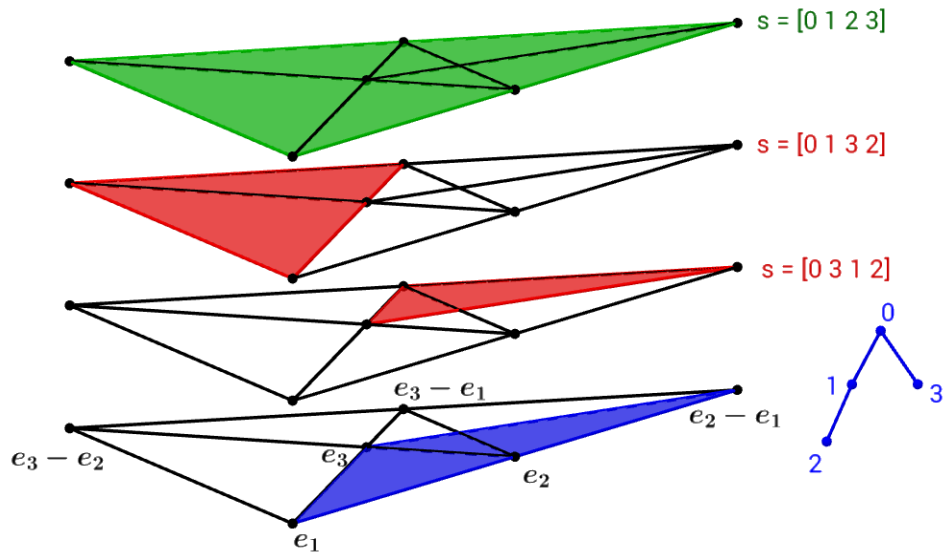

**Figure 3**: A representation by three path cones

However, for the tree $T$ shown in Figure 2(b), the relevant permutations are then [0213] and [0231], both with negative distortion. So, $\chi_T$ is now the *sum* of the two path cones, as shown in Figure 2(b).

Of the remaining seven trees, six are of the form shown in Figure 3. Now, there are three completions of $\prec_T$ to a full order. For example, for the tree $T$ shown in Figure 3, the three paths that are compatible with $T$ correspond to $s_1 = [0123]$, $s_2 = [0132]$ and $s_3 = [0312]$. Again, $\mathbf{d}(T) = 0$, while $\mathbf{d}(T_{s_i})$ is $0, 1, 1$, for $i = 1, 2, 3$, respectively. Figure 3 shows then the decomposition, with the original $\chi_T$ colored blue, the positive cone colored green and the negative ones colored red.

Finally, the tree $T = (e_i : i \in [1{:}3])$ shown in Figure 4 is compatible with all the six path trees, and has $\mathbf{d}(T) = 0$. Its decomposition is shown in Figure 4, with green representing the cones with positive coefficients and red those with negative coefficients. The result itself, i.e., the positive octant, is shown in blue.

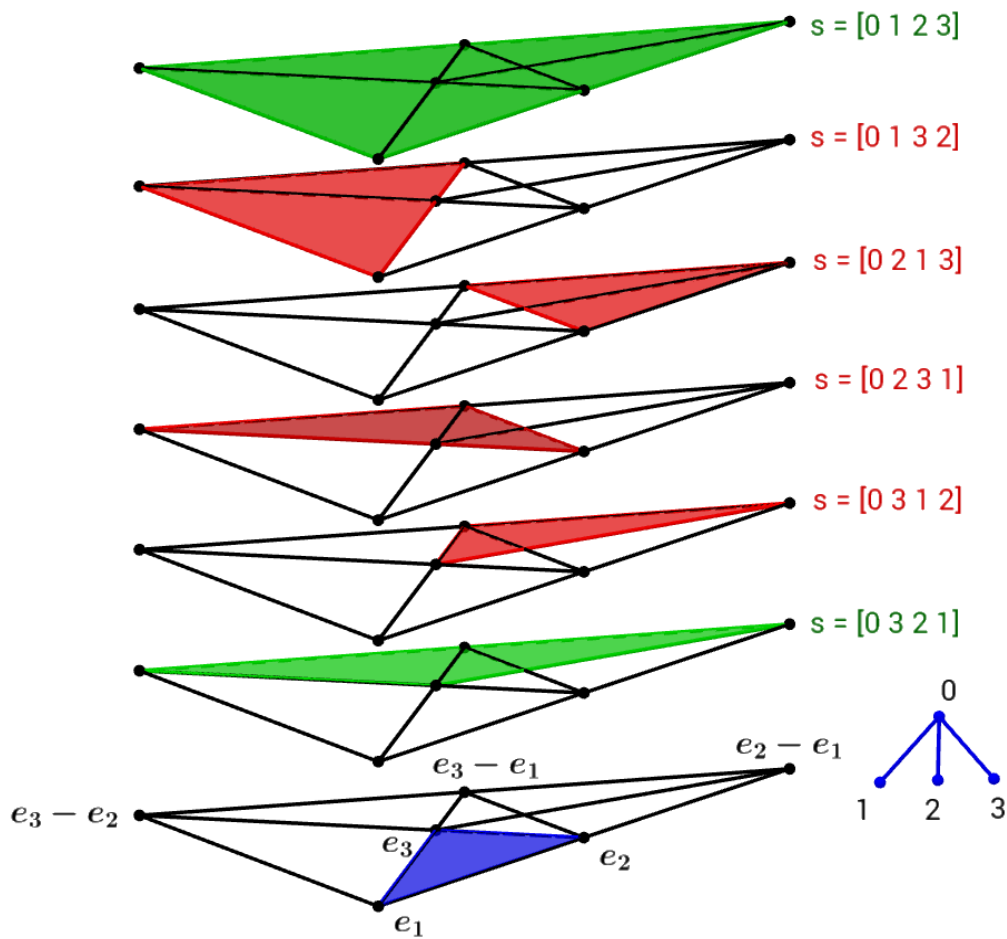

Figure 4: The representation of the positive octant

## 3. Proofs

We work with polynomials in the $n$ variables
$$t = (t(1), \ldots, t(n)).$$

An edge $x = e_j - e_i \in G$ is then associated with the linear homogeneous polynomial

$$p_x(t) = t(j) - t(i).$$

Thus, the differential operator

$$p_x(D)$$

is the directional derivative in the $x$ direction. If $Y \subset G$, i.e., a subset of the edge set of $G$, oriented as in $G$, we denote

$$p_Y := \prod_{x \in Y} p_x.$$

Thus, for every $T \in \mathbb{B}(G)$, $p_{G \backslash T}$ is a homogeneous polynomial of degree $\binom{n+1}{2} - n = \binom{n}{2}$. Denote

$$\mathrm{soc}(\mathcal{P}(G)) := \mathrm{span}\{p_{G \backslash T} : T \in \mathbb{B}(G)\}.$$

We have then

**Result 3.1.**
(1) $\dim \mathrm{soc}(\mathcal{P}(G)) = n!$.
(2) For every path tree $T_s$, $s \in S_n$, there exists a polynomial $\mathbf{M}_s$, homogeneous of degree $\binom{n}{2}$, such that:
(2a) If $T \in \mathbb{B}(G)$, and $Y := G \backslash T$, then $p_{Y^s}(D)\mathbf{M}_s \in \{0, 1\}$.
(2b) For any $Y \subset G$, $p_Y(D)\mathbf{M}_s = 0$ if and only if the graph $G^s \backslash Y^s$ has a source other than 0.

**Proof:** (1) above is well known. It is obtained by combining Theorem 7.2 in [6] with some basics of algebraic graph theory, [4]. See also [5], [3], [1] and [7].

(2) follows from the results of [7] which are reviewed in Appendix A. □

**Corollary 3.2.** *The polynomials*

$$\mathbf{P}_s := p_{G^s \backslash T_{s,on}}, \quad s \in S_n,$$

*form a basis for* $\mathrm{soc}(\mathcal{P}(G))$. *The polynomials* $(\mathbf{M}_s : s \in S_n)$ *are biorthogonal to that basis via the pairing*

$$\Pi \ni (p, q) \mapsto \langle p, q \rangle := p(D)q(0).$$

**Comment.** The basis $(\mathbf{P}_s)$ is not new. It appears in [3], and is utilized in [2] and [1]. However, we will need the special properties of the dual basis $(\mathbf{M}_s)$, hence provide a complete proof, starting, indeed, with the dual basis.

**Proof.** Fix $s, s' \in S_n$. Then, by (2b) of Result 3.1, $\mathbf{P}_{s'}(D)\mathbf{M}_s = 0$ iff $T^s_{s'}$ has more than one source. But, for any tree $T$, only the natural orientation $T_{on}$ is void of a source $i \in [1:n]$. So, $\mathbf{P}_{s'}(D)\mathbf{M}_s = 0$, unless $T^s_{s'} = T_{s',on}$. This is possible iff $s = s'$. So, $\mathbf{P}_{s'}(D)\mathbf{M}_s = 0$ iff $s \neq s'$. In the complementary case, (2a) of Result 3.1 implies that $\mathbf{P}_s(D)\mathbf{M}_s \in \{0, 1\}$, and since we argued already that it is not 0 then it is 1.

Thus the biorthogonality relation is proved, and it then implies that $(\mathbf{P}_s : s \in S_n)$ are linearly independent. Since there are $n!$ of them, we conclude from (1) of Result 3.1 that they form a basis for $\mathrm{soc}(\mathcal{P}(G))$. □

In the theorem below, we denote, for any $s \in S_n$ and for any $Y \subset G$, by

$$\mathrm{mm}(Y, s),$$

the mismatch between $Y$ and $Y^s$, i.e., the number of edges whose orientations in these two oriented sets differ.

**Theorem 3.3.** *Let $T \in \mathbb{B}(G)$. Then*

$$p_{G\setminus T} = \sum_{s \in S_n} (-)^{\mathrm{mm}(G\setminus T, s)} C(T, T_s)\, \mathbf{P}_s.$$

**Proof:** By Corollary 3.2,

$$p_{G\setminus T} = \sum_{s \in S_n} \langle p_{G\setminus T}, \mathbf{M}_s \rangle\, \mathbf{P}_s.$$

Fixing $s \in S_n$, we want first to determine whether the coefficient of $\mathbf{P}_s$ is 0 or not. To this end, denote $Y := G\setminus T$, and apply (2b) of Result 3.1. Since $G^s \setminus Y^s = T^s$, the condition there is that $T^s$ does not have a source other than 0, which is equivalent to $T^s = T_{on}$. So, $\langle p_{G\setminus T}, \mathbf{M}_s \rangle \neq 0$ iff $T^s = T_{on}$, which is equivalent to the compatibility of $T$ and $T_s$, i.e, to the condition $C(T, T_s) = 1$.

Now assume that, indeed, $C(T, T_s) = 1$, hence that $T^s = T_{on}$. We have just argued that in that case $p_{Y^s}(D)\mathbf{M}_s \neq 0$, hence by (2a) of Result 3.1, $p_{Y^s}(D)\mathbf{M}_s = 1$. Since we apply $p_Y(D)$ instead of $p_{Y^s}(D)$, we need to flip the sign of the result, once for each mismatch count in $\mathrm{mm}(Y, s)$. The result thus follows. □

We now need to transport all the above results from the space $\mathrm{soc}(\mathcal{P}(G))$ to the simplicial space $\mathcal{T}(G)$. To this end, we utilize the truncated power associated with $G$, $\mathrm{TP}_G$, [5]. It is a piecewise-polynomial function supported in the cone $\mathrm{pos}(G)$. For our purposes, the only property of $\mathrm{TP}_G$ that we need is the following

**Result 3.4.** *Let $T \in \mathbb{B}(G)$, then*

$$p_{G\setminus T}(D)\mathrm{TP}_G = \chi_T.$$

*So, the space $\mathcal{T}(G)$ is isomorphic to the space $\mathrm{soc}(\mathcal{P}(G))$ via the map*

$$\mathrm{soc}(\mathcal{P}(G)) \ni p \mapsto p(D)\mathrm{TP}_G \in \mathcal{T}(G).$$

*Thus, Proposition 1.1 follows from (1) of Result 3.1, while Theorem 1.2 follows from Corollary 3.2. Finally, Theorem 1.4 follows from Theorem 3.3 as we argue now:*

**Proof of Theorem 1.4.**

Fix $T \in \mathbb{B}(G)$. Then, by Result 3.4, when combined with Corollary 3.2,

$$\chi_T = p_{G\setminus T}(D)\mathrm{TP}_G = \sum_{s \in S_n} \langle p_{G\setminus T}, \mathbf{M}_s \rangle\, \mathbf{P}_s(D)\mathrm{TP}_G.$$

Now, Theorem 3.3 provides the coefficient exactly:

$$\langle p_{G\backslash T}, \mathbf{M}_s \rangle = (-)^{\mathrm{mm}(G\backslash T, s)} C(T, T_s). \tag{3.5}$$

However, Result 3.4 asserts that $\chi_s = p_{G\backslash T_s}(D)\mathrm{TP}_G$, while $\mathbf{P}_s = p_{G^s \backslash T_s^s}$, so we need to account for the sign change when passing from the latter operator to the former, i.e., we need to account for the parity of the number

$$\mathrm{mm}(G\backslash T_s, s) = \mathrm{mm}(G, s) - \mathrm{mm}(T_s, s) = \mathrm{mm}(G, s) - \mathbf{d}(T_s).$$

The number we already have in (3.5) is

$$\mathrm{mm}(G\backslash T, s) = \mathrm{mm}(G, s) - \mathrm{mm}(T, s).$$

So, the parity of the set

$$\mathbf{d}(T_s) + \mathrm{mm}(T, s)$$

determines the sign of the coefficient, whenever such sign is required, i.e., whenever $C(T, T_s) = 1$. But this latter condition is equivalent to $T^s = T_{on}$, and in that case, $\mathrm{mm}(T, s) = \mathbf{d}(T)$. □

## 4. Appendix A: on the construction of the dual polynomials $\mathbf{M}_s$

We outline the construction from [7], save the explicit algorithm for the construction of these polynomials: for the sake of this paper, we need the *existence* of such polynomials, with properties as in Result 3.1, and not their explicit form.

One starts with a graph $G$, this time undirected, and connected. The dual polynomials are indexed then (in lieu the indexing by $S_n$) by the set

$$\mathrm{O}(G)$$

of all the acyclic orientations of $G$ with a single source 0. The polynomials

$$\mathcal{M} := \{\mathbf{M}_{G_\iota} : G_\iota \in \mathrm{O}(G)\},$$

are all homogeneous of degree $\#G - n$, with $G$ identified with its edge set. Fixing $G_\iota \in \mathrm{O}(G)$, one partitions the edge set $G$ into $n$ subsets $X_{G_\iota, i}$, $i \in [1{:}n]$, as follows:

$$X_{G_\iota, i} := \{x \in G : \ x \text{ is an } \textit{inflow} \text{ edge of } i \text{ in } G_\iota\}.$$

So, if $x = e_i - e_j \in G_\iota$ it goes to $X_{G,i}$. [7] proves that, up to normalization, there exists a unique polynomial of degree $\#G - n$ in the joint kernel of the differential operators

$$p_{X_{G,i}}(D), \quad i \in [1{:}n].$$

Then $\mathbf{M}_{G_\iota}$ is defined this way up to normalization, and we refer to [7] for the normalization details.

Now, if $x = e_i - e_j \in G_\iota$ and $G_\iota' := G_\iota \backslash x$, then it follows easily that: (1) if $X_{G,i} = \{x\}$ then $p_x(D)\mathbf{M}_{G_\iota} = 0$, and the graph $G_\iota'$ has then $i$ as an additional source. (2) Otherwise, the graph $G_\iota'$ has still only one source, 0, and the new differential operators

$$p_{X_{G_\iota',i}}(D), \quad i \in [1{:}n],$$

annihilate, each, $p_x(D)\mathbf{M}_{G_\iota}$. Then, once one proves that $p_x(D)\mathbf{M}_{G_\iota} \neq 0$, the uniqueness assertion above implies that, up to normalization

$$\mathbf{M}_{G_\iota'} = D_x \mathbf{M}_{G_\iota}.$$

Property (2) in Result 3.1 follows then from a slightly more careful argument, when appropriate normalizations in the definition of $\mathbf{M}_{G_\iota}$ are applied.

Note that the definition given above for the dual polynomials easily implies that, once $G$ is itself a spanning tree $T$, $\mathrm{O}(G) = \{T_{on}\}$, and $\mathbf{M}_{T_{on}} = 1$ (up to normalization).